\documentclass[11pt, a4paper]{amsart}

\usepackage[left=3cm, top=3cm,bottom=3cm,right=3cm]{geometry}
\usepackage{graphicx}

\usepackage[utf8]{inputenc}
\usepackage{lmodern}
\usepackage{amssymb,amsfonts,amsthm,amsmath,bbm,mathrsfs,aliascnt,mathtools}
\usepackage[english]{babel}
\usepackage[colorlinks]{hyperref}

\usepackage{ulem}
\usepackage{nicefrac}

\theoremstyle{plain}
\newtheorem*{thm}{Theorem}

\newtheorem{lem}{Lemma}

\DeclareMathOperator{\Diff}{Diff}

\begin{document}

\title{Distortion in the group of circle homeomorphisms}
\author{Juliusz Banecki}
\address{Juliusz Banecki, Institute of Mathematics Polish Academy of Sciences, Abrahama 18, 81-967 Sopot, Poland}
\email{juliusz.banecki@autonomik.pl}

\author{Tomasz Szarek}
\address{Tomasz Szarek, Faculty of Physics and Applied Mathematics, Gda{\'n}sk University of Technology,
ul. Gabriela Narutowicza 11/12, 80-233 Gda{\'n}sk \& Institute of Mathematics Polish Academy of Sciences, Abrahama 18, 81-967 Sopot, Poland}
\email{tszarek@impan.pl}

\maketitle

\begin{abstract}
Let $G$ be the group ${\rm PAff}_+({\bf S}^1)$ of piecewise--affine circle homeomorphisms or the group  ${\Diff}^{\infty}(\mathbb R/\mathbb Z)$ of smooth circle diffeomorphisms.  A constructive proof that all irrational rotations are distorted in $G$ is given. 
\end{abstract}


\section{Introduction}

Let $G$ be a group with some finite generating set $\mathcal G$. We define the metric $d_{\mathcal G}$ on $G$ by taking $d_{\mathcal G} (g_1, g_2)$ to be the infimum over all
$k\ge 0$ such that there exist $f_1,\ldots, f_k\in\mathcal G$ and $\epsilon_1,\ldots,\epsilon_k\in\{-1, 1\}$ satisfying $g_2=f_1^{\epsilon_1}\cdots f_k^{\epsilon_k}g_1$.

Now let $H$ be an arbitrary group. An element $f\in H$ is called {\it distorted} in $H$  if there exists
a finitely generated subgroup $G\subset H$ containing $f$ such that
$$
\lim_{n\to\infty}\frac{d_{\mathcal G}(f^n, {\rm id})}{n}=0
$$ 
for some (and hence for every) generating set $\mathcal G$. Since the limit always exists, it is enough to verify it for some subsequence. 
The notion of distortion comes from geometric group theory and it was introduced by M. Gromov in \cite{Gromov}.
\vskip2mm

The problem of the existence of distorted elements in some groups of homeomorphisms has been intensively studied for many years (see \cite{C_F, Franks, F_H, G_L, Navas1, Pol}). Substantial progress has been achieved for groups of diffeomorphisms of manifolds. In particular, A. Avila \cite{Avila} proved that rotations with irrational rotation number are distorted in the group of smooth diffeomorphisms of the circle.  In this note we give a constructive proof that all irrational rotations are distorted both in the group of piecewise-affine circle homeomorphisms, ${\rm PAff}_+({\bf S}^1)$, and in the group of smooth circle diffeomorphisms, ${\Diff}^{\infty}(\mathbb R/\mathbb Z)$. 
The result gives an answer to {\bf Question 11} in \cite{Navas} (see also {\bf Question 2.5} in \cite{Franks}). So far it has not even been known whether there exist distorted elements in ${\rm PAff}_+({\bf S}^1)$. Now from \cite{G_L} it follows that distorted elements, if they exist, are conjugate to rotations.

\vskip2mm
From now on let $G$ be either the group ${\rm PAff}_+({\bf S}^1)$ or 
 ${\Diff}^{\infty}(\mathbb R/\mathbb Z)$. We say that $g\in G$ is {\it trivial on some set} if there exists a non-empty open set $I\subset {\bf S}^1$  such that $g(x)=x$ for $x\in I$. The set of all homemomorphisms in $G$ which are trivial on some set will be denoted by $G_{{\rm {triv}}}$. By ${\rm T}$ we denote the set of all rotations.
 
 \vskip2mm
This note is devoted to the proof of 

\begin{thm}
All irrational rotations are distorted in $G$.
\end{thm}

\section{Proofs}

We first formulate two lemmas and deduce the theorem. 
The proofs of the lemmas will be given at the end of the paper.

\begin{lem}\label{L1}
For any irrational  rotation $T_{\alpha}$ and $g\in G_{{\rm {triv}}}\cup {\rm T}$  there exist a finite generating set $\mathcal G_g\subset G$ and a constant $C>0$ such that 
\[
d_{\mathcal G}(T_{\alpha}^n g T_{\alpha}^{-n}, {\rm id})\le C\log n\quad\text{for all $n\ge 1$}.
\]
\end{lem}

\vskip2mm

\begin{lem}\label{Prop1_18.04.21}  In $G$ there exist
$g_1,\ldots, g_l\in G_{{\rm {triv}}}\cup{\rm T}$ and $k, k_1,\ldots, k_l\in {\mathbb Z}$
with $k\neq k_1+\cdots+k_l$, such that for each sufficiently small $\beta>0$ the element
$x=T_{\beta}$ satisfies 
\begin{equation}\label{e1_17.04.21}
x^{k_1}g_1x^{k_2}g_2\cdots x^{k_l}g_l=x^k.
\end{equation} 
\end{lem}

{\it Proof of the theorem.}  Fix an irrational rotation $T_{\alpha}$. From Lemma \ref{Prop1_18.04.21} it follows that in $G$ there exists an equation of the form (\ref{e1_17.04.21}) such that $x=T_{\beta}$, for all sufficiently small $\beta$, is its solution. Let $\mathcal G=\mathcal G_{g_1}\cup\cdots\cup \mathcal G_{g_l}$, where $\mathcal G_{g_i}$, $i=1,\ldots, l$, are finite generating sets derived from Lemma \ref{L1} for $T_{\alpha}$.
We may rewrite equation (\ref{e1_17.04.21}) in the form
\begin{equation}\label{e2_18.04.21}
x^{k_1}g_1x^{-k_1} x^{k_2+k_1}g_2x^{-k_2-k_1}\cdots x^{k_1+\cdots+k_l}g_l x^{-k_1-\cdots-k_l}=x^{k-k_1-\cdots-k_l}.
\end{equation}
Let $\beta_0$ be a positive constant such that $x=T_{\beta}$ for $\beta\in (0, \beta_0)$ satisfies (\ref{e2_18.04.21}). Set $m:=k-k_1-\cdots-k_l$, and let $(n_i)$ be an increasing sequence of  integers such that $n_i\alpha \in (0, \beta_0)\, ({\rm mod}\, 1)$. From Lemma \ref{L1} it follows that
$$
d_{\mathcal G}(T_{\alpha}^{n_i(k_1+\cdots+k_j)} g_j T_{\alpha}^{-n_i(k_1+\cdots+k_j)}, {\rm id})\le C_j\log n_i\quad\text{for all $i\ge 1$ and $j=1,\ldots, l$}.
$$
Since $x=T_{n_i\alpha}$ satisfies (\ref{e2_18.04.21}), we obtain
$$
d_{\mathcal G}(T_{\alpha}^{n_i m}, {\rm id})\le \sum_{j=1}^l C_j\log n_i:= C \log n_i\quad\text{for all $i\ge 1$.}
$$
Hence 
$$
\lim_{n\to\infty}\frac{d_{\mathcal G}(T_{\alpha}^n, {\rm id})}{n}=
\lim_{i\to\infty}\frac{d_{\mathcal G}(T_{\alpha}^{n_i m}, {\rm id})}{m n_i}\le \lim_{i\to\infty}\frac{C \log n_i}{m n_i}=0
$$ 
and the proof is complete. $\square$
\vskip4mm

{\it Proof of Lemma 1.} The proof relies on the observation that for a given interval $I\subset (0, 1)$  there exists a finite generating set $\mathcal G\subset G$ such that for any $n\ge 1$ there exists a 
homemorphism $h_n$ with $d_{\mathcal G}(h_n, {\rm id})\le C\log n$ for some constant $C>0$ independent of $n$, and $h_n(x)=T_{\alpha}^n(x)$ for $x\notin I$. 
Without loss of generality we may assume that $I=(a, 1).$ Let $m\ge 1$ be an integer such that $a+2/m<1$. Let $h\in G $ be any homeomorphism such that $h(x)=x/2$ for $x\in [0, a+2/m)$, and let $r(x)=x+1/m$. 

We shall define $h_n$ by induction. Set $h_0=h$. If $n$ is odd we put $h_n=T_{\alpha} h_{n-1}$. If $n$ is even, we take $s_n:=h_{n/2} h$ 
and observe that $s_n((0, a))=(n\alpha/2, a/2+n\alpha/2)$. Let $k\in\{1,\ldots, m\}$ be such that $n\alpha/2+k/m\in [0, 1/m)$\,(mod $1$). Then $r^k s_n((0, a))\subset (0, a/2+1/m)$. Therefore 
$$
h^{-1} r^k h_{n/2} h(x)=2(x/2+n\alpha/2+k/m)=x+n\alpha+2k/m=T_{\alpha}^n(x)+2k/m
$$ 
for $x\in (0, a)$. Put $h_n:=r^{-2k}h^{-1} r^k h_{n/2} h$, and let $\mathcal G:=\{T_{\alpha}, h, r\}$. Note that
$$
d_{\mathcal G}(h_n, {\rm id})\le 3m+3+d_{\mathcal G}(h_{\lfloor n/2\rfloor}, {\rm id}).
$$
Thus we obtain $d_{\mathcal G}(h_n, {\rm id})\le C\log n$. Finally, observe that
for any $g\in G_{{\rm {triv}}}$ such that $g(x)={\rm id}$ on $I$ we have
$$
T_{\alpha}^n g T_{\alpha}^{-n}=h_n g h_n^{-1}.
$$
Therefore, we obtain
$$
d_{\mathcal G}(T_{\alpha}^n g T_{\alpha}^{-n}, {\rm id})\le C\log n.
$$
In the case when $g$ is a rotation the conclusion of the lemma is obvious. $\square$
\vskip4mm
{\it Proof of Lemma 2.} Let $\beta\in (0, 10^{-3})$, and let $f_1\in G_{{\rm {triv}}}$ be arbitrary such that
$$
f_1(x)=0.4+2(x-0.4)\,\,\text{for $x\in [0.4, 0.6]$ and } f_1(x)=x\,\,\text{for $x\in [0.9, 1.1]$}.
$$
Set
$$
H_1=T_{2\beta}^{-1} f_1 T_{2\beta} f_1^{-1}.
$$
It is obvious that
$$
H_1(x)=x+2\beta\,\,\text{for $x\in [0.41, 0.79]$  and }  H_1(x)=x\,\,\text{for $x\in [0.91, 1.09]$}.
$$
Define 
$$
H_2=T_{1/2} H_1^{-1} T_{1/2} H_1,
$$
and observe that 
$$
H_2(x)=x-2\beta\quad\text{for $x\in [0.95, 1]$}.
$$
Simple computation gives 
$$
T_{1/2} H_2T_{1/2} H_2={\rm id}.
$$
Set 
$$
H_3 =T_{2\beta} H_2.
$$
Then we have 
$$
H_3(x)=x\quad\text{for $x\in [0.95, 1]$}
$$
and
\begin{equation}\label{e1_22.04.21}
T_{2\beta+1/2}H_3 T_{-2\beta-1/2} H_3=T_{4\beta}.
\end{equation}
Take an arbitrary $f_2\in G_{{\rm {triv}}}$ satisfying 
$$
f_2(x)=2x\quad\text{for $x\in [0, 0.49]$},
$$
and define
$$
H_4=f_2^{-1} H_3f_2.
$$
It is easy to see that
$$
H_4(x)=
\begin{cases}
H_3(2x)/2 & \text{for }  x\in [0, 1/2),\\
x & \text{for }   x\in [1/2, 1).
\end{cases}
$$
Let
\begin{equation}\label{e2_22.04.21}
H_5=T_{1/2} H_4 T_{1/2}H_4.
\end{equation}
Observe that the graph of $H_5$ is built 
from two scaled copies of $H_3$, i.e.
$$
H_5(x)=
\begin{cases}
H_3(2x)/2 & \text{for }  x\in [0, 1/2),\\
H_3(2x-1)/2+1/2 & \text{for }  x\in [1/2, 1).
\end{cases}
$$
Therefore, by (\ref{e1_22.04.21}) and (\ref{e2_22.04.21}), we finally obtain
\begin{equation}\label{e1_18.04.21}
T_{\beta+1/4}H_5 T_{-\beta-1/4}H_5=T_{2\beta}. 
\end{equation}
Indeed, this is easy to see if we realize that (\ref{e1_18.04.21}) is simply equation (\ref{e1_22.04.21}) rewritten in the new coordinates $(x/2, y/2)$.
Plugging subsequently $H_5, H_4, H_3, H_2$ and $H_1$ into formula (\ref{e1_18.04.21}) we obtain the required equation with $g_1,\ldots, g_l\in\{f_1, f_2, f_1^{-1}, f_2^{-1}, T_{1/2}, T_{-1/2}, T_{1/4}, T_{-1/4}\}\subset G_{{\rm {triv}}}\cup {\rm T}$. Moreover, on the left hand side of this equation
there are eight more $T_{\beta}$'s than $T^{-1}_{\beta}$'s. On the other hand, on the right hand side we have  two $T_{\beta}$'s. Since $\beta\in (0, 10^{-3})$ was arbitrary and the functions $g_1,\ldots, g_l$ were independent of $\beta$, the proof of the lemma is complete. $\square$

\end{document}